\documentclass{amsart}
\usepackage{graphicx, latexsym, amsmath, amssymb}
\newtheorem{thm}{Theorem}[section]

\newtheorem{lem}[thm]{Lemma}
\newtheorem{prop}[thm]{Proposition}
\theoremstyle{definition}
\newtheorem{defn}[thm]{Definition}
\theoremstyle{remark}
\newtheorem{rem}{Remark}
\newtheorem{exmp}{Example}
\numberwithin{equation}{section}
\newcommand{\scalar}[2]{\langle#1,#2\rangle}
\newcommand{\Real}{\mathbb{R}}
\newcommand{\Complex}{\mathbb{C}}
\newcommand{\To}{\longrightarrow}

\newcommand{\pars}[1]{\left(#1\right)}
\newcommand{\brck}[2]{\{#1,#2\}}
\newcommand{\set}[1]{\left\{#1\right\}}
\newcommand{\jami}[2]{\sum\limits_{#1}^{#2}}
\newcommand{\all}[1]{\forall\,#1}
\newcommand{\vbar}{\Bigr|}
\newcommand{\gamomd}{\,\Rightarrow\,}

\newcommand{\smooth}[1]{\ensuremath{C^\infty(#1)}}
\newcommand{\fields}[2]{\mathcal{V}^{#1}(#2)}

\newcommand{\tolfas}{\,\Leftrightarrow\,}
\newcommand{\lder}[2]{L_{#1}(#2)}
\newcommand{\talg}[1]{\widetilde{#1}}
\newenvironment{prf}{\textbf{Proof.}}{$\blacksquare$}
\newcommand{\imply}{\Rightarrow}
\begin{document}
\title[]{Algebraic definition of Holonomy on Poisson Manifold}%
\author{Zakaria Giunashvili}%
\address{Department of Theoretical Physics, A. Razmadze Institute of Mathematics, 1 Alexidze, Tbilisi ge-0193, Georgia}%
\email{zaqro@mes.gov.ge}%

\thanks{The work was supported by the Georgian National Science Foundation under grant
GNSF/ST-06/4-050}%

\date{June 3, 2007}%
% ----------------------------------------------------------------
\begin{abstract}
We give an algebraic construction of connection on the symplectic leaves of
Poisson manifold, introduced in \cite{Ginzburg}. This construction is suitable
for the definition of the linearized holonomy on a regular symplectic
foliation.
\end{abstract}
\maketitle
% ----------------------------------------------------------------
%%%%%%%%%%%%%%%%%%%%%%%%%%%%%%%%%%%%%%%%%%%%%%%%%%%
\section{Introduction}
If the rank of a Poisson structure on a smooth manifold $M$ is constant, there
are various tools in the theory of foliations that can be applied to study the
corresponding symplectic foliation. But there are many even very simple
degenerate Poisson structures that are singular, i.e, the rank of the
corresponding contravariant tensor field is not constant. Therefore, some new
methods are needed to handle this situation.

In this paper we give an algebraic coonstruction of the linear holonomy on a
Poisson manifold. This construction is suitable for the case, when the Poisson
structure is not regular and the necessary requirements for the Poisson
structure, directly follow from the construction.
%%%%%%%%%%%%%%%%%%%%%%%%%%%%%%%%%%%%%%%%%%%%%%%%%%%
\section{Algebraic Definition of Bott Connection}
\newcommand{\der}{Der(A)}
\newcommand{\deri}{Der_I(A)}
\newcommand{\derio}{Der_I(A)_0}
\newcommand{\gami}{\Gamma(I)}
\newcommand{\gamio}{\Gamma(I)_0}
\newcommand{\inter}{\cap}
%%%%%%%%%%%%%%%%%%%%%%%%%%%%%%%%
Let $A$ be a commutative algebra and $Der(A)$ be the Lie algebra of derivations
of $A$. As $A$ is a commutative algebra, the space $Der(A)$ is a module over
$A$. For a given ideal $I$ in $A$ let us define the following Lie subalgebras
in $Der(A)$:
$$
\deri=\set{X\in\der\ |\ X(I)\subset I}
$$
and its subalgebra
$$
\derio=\set{Y\in\der\ |\ Y(A)\subset I}.
$$
It is clear that both of them are $A$-submodules in $Der(A)$.

There is a canonical homomorphism
$$
\rho:\deri\To Der(A/I)
$$
and as it follows from the definition, we have that $\derio=\ker(\rho)$. When
$\rho$ is an epimorphism, the ideal $I$ is said to be a \emph{submanifold
constraint ideal} (see \cite{Masson}). For brevity we shall use the term
\emph{submanifol ideal}. In this case we have the isomorphism
$$
Der(A/I)\cong\frac{\deri}{\derio}
$$
We use the term \emph{distribution} for any Lie subalgebra which is also an
$A$-submodule of $\der$.
\begin{defn}[Integral of Distribution]
An ideal $I$ in $A$ is said to be an \emph{integral} for a given distribution
$D$ if $D\subset\deri$.
\end{defn}
Let us denote the intersection $D\inter\derio$ by $D_0$.
\begin{defn}[Complete Integral]
An ideal $I\subset A$ is said to be a \emph{complete integral} for a given
distribution $D\subset\der$ if the following three conditions are satisfied
\begin{itemize}
\item[(1)] $I$ is a submanifold ideal for $A$
\item[(2)] $I$ is an integral for $D$
\item[(3)] $D/D_0=\frac{\deri}{\derio}\quad$ (\ $\cong Der(A/I)$ because of (1))
\end{itemize}
\end{defn}
Let $I\cdot D=\set{\sum a_iX_i\ |\ a_i\in I,\ X_i\in D}$. It is clear that
$I\cdot D$ is a subspace of $D_0$.
\begin{defn}[Regular Integral]
A complete integral $I$ for a given distribution $D\in\der$ is said to be
\emph{regular} if $D_0=I\cdot D$.
\end{defn}
The above definition of a regular integral is in agreement with the following
theorem from the classical theory of distributions
\begin{thm}[see \cite{Dazord}]
Let $D$ be a completely integrable distribution on a smooth manifold $M$, and
$x_0\in M$ be a point where $\dim(D_{x_0})=r$. Then the point $x_0$ has an open
neighborhood $U$ with coordinates $(u_1,\ldots,u_n)$ such that
$$
\set{\frac{\partial}{\partial u_i}(x)\ |\ i=1,\ldots,r}\subset D_x,\ \all x\in
U
$$
\end{thm}
It follows from this theorem that if $x_0$ is a regular point for $D$, then
$\dim(D_{x_0})$ is locally maximal, and therefore, the set
$\set{\frac{\partial}{\partial u_i}\ |\ i=1,\ldots,r}$ is a local basis of the
$\smooth{M}$-module $D$. If $N\subset M$ is the integral submanifold for $D$,
passing through the point $x_0$ and $X\in D$ is such that $X|_N=0$, then
$X=\jami{i=1}{r}\varphi_i\frac{\partial}{\partial u_i},\ \varphi_i\in I_N$,
where $I_N$ denotes the ideal of such smooth functions on $M$ that vanish on
$N$.

For a given ideal $I\in A$ let $\pi:A\To A/I$ be the natural quotient map.
Introduce the following $A$-modules
$$
\gami=\set{X:A\To A/I\ |\ X(ab)=X(a)\pi(b)+\pi(a)X(b)}
$$
and its submodule
$$
\gamio=\set{Y\in\gami\ |\ Y(I)=\set{0}}
$$
We call the quotient module
$$
V(I)=\frac{\gami}{\gamio}
$$
the space of \emph{transversal derivations}. It is a module over the quotient
algebra $A/I$. Consider the following operation $[\ ,\
]:\deri\times\gami\To\gami$, defined as
$$
[U,s]=\rho(U)\circ s-s\circ U,\quad U\in\deri,\ s\in\gami
$$
The submodule $\gamio$ is invariant for this operation and therefore we can
reduce this operation to
$$
[\ ,\ ]:\deri\times V(I)\To V(I)
$$
The following properties of this operation follow directly from the definition
\begin{equation}\label{bracket_props}
\begin{array}{ll}
(1) & [aU,v]=a[U,v]\\
& \\
(2) & [U,kv]=U(k)v+k[U,v]
\end{array}
\end{equation}
where $a\in A,\ k\in A/I,\ U\in\deri$ and $v\in V(I)$.
\begin{defn}[Linear Connection]
\emph{Linear connection} for $I\subset A$ which is a \emph{regular integral}
for a distribution $D\in\der$, is an action (covariant derivation) of the
elements of Lie algebra $Der(A/I)\cong\frac{D}{I\cdot D}$ on the $A/I$-module
$V(I)$, defined as
$$
\nabla_X(v)=[\talg{X},v]
$$
where $X\in\frac{D}{I\cdot D},\ v\in V(I)$ and $\talg{X}\in D$ is an extension
of $X$: $\rho(\talg{X})=X$.
\end{defn}
\textbf{Independence from the extension $\talg{X}$:} if $\talg{X}_1$ and
$\talg{X}_2$ in $D$ are such that $\rho(\talg{X}_1)=\rho(\talg{X}_2)$ then we
have
$$
\begin{array}{c}
\talg{X}_1-\talg{X}_2\ \in I\cdot D\tolfas\talg{X}_1-\talg{X}_2=\sum k_iU_i,\ k_i\in I,\ U_i\in D\gamomd\\
\\
\gamomd[\talg{X}_1-\talg{X}_2,v]=[\sum k_iU_i,v]=\sum k_i[U_i,v]=0\ \in\
V(I)=\frac{\gami}{\gamio}
\end{array}
$$
\medskip\\
The properties of covariant derivation for $\nabla$ easily follow from the
properies of the operation $[\ ,\ ]$. In the case of a singular (i.e., not
regular) integral, the definition of linear connection encounters the following
problem: the subalgebra $D_0$, in general, does not coincide with $I\cdot D$,
i.e., not every $\talg{X}\in D$ that vanishes on $A/I$ is of the form
$\talg{X}=\sum k_iU_i$ with $k_i\in I$ and $U_i\in D$. Hence the independence
of the expression $[\talg{X},v]$ from the extension $\talg{X}$ is problematic.
A Poisson structure helps to somehow resolve the difficulty. Some authors (see
\cite{Ginzburg}) define the covariant derivation of a transversal vector field
not by tangent vectors but by 1-forms. In the next section we shall review this
construction in the case of a Poisson manifold.
%%%%%%%%%%%%%%%%%%%%%%%%%%%%%%%%%%%%%%%%%%%%%%%%%%%
\section{Linear Connection on a Symplectiv Leaf of Poisson Manifold}
\newcommand{\ortho}{T(N)^\perp}
\newcommand{\restr}{T^*(M)|_N}
%%%%%%%%%%%%%%%%%%%
Let $(M,\brck{\ }{\ })$ be a Poisson manifold. For any function
$f\in\smooth{M}$, we denote by $X_f$ the Hamiltonian vector field on $M$
corresponding to $f$. The antisymmetric contravariant tensor (bivector) field
on $M$ corresponding to the Poisson structure, defines a homomorphism of the
vector bundles $T^*(M)\To T(M)$. We denote the tangent vector corresponding to
the cotangent vector $\alpha\in T^*(M)$ also by $X_\alpha$. For a differential
1-form $\omega$ we have a vector field $X_\omega$ defined by the above
homomorphism.

Let $N\subset M$ be a symplectic leaf. Denote by $\ortho$ the subbundle of the
restricted bundle $\restr$:
$$
\ortho=\set{\omega\in T_x^*(M)\ |\ x\in N,\ \omega(T_x(N))=0}
$$
\emph{Linear connection} on the vector bundle $\ortho$ is a covariant
derivation of the sections of this bundle by the sections of the vector bundle
$\restr$ defined as follows
$$
\nabla_\alpha(\omega)=\lder{X_{\talg{\alpha}}}{\talg{\omega}}|_N
$$
where:
\begin{itemize}
\item
    $\alpha$ is a section of $\restr$ and $\talg{\alpha}$ is its extension (at least local) on $M$:
    $\talg{\alpha}|_N=\alpha$;
\item
    $\omega$ is a section of $\ortho$ and $\talg{\omega}$ is its extension (at least local) on M:
    $\talg{\omega}|_N=\omega$;
\item
    $\lder{\cdot}{\cdot}$ is the operation of Lie derivation.
\end{itemize}
If $\talg{\alpha}=\varphi\cdot\talg{\alpha}_1$, where $\varphi\in\smooth{M}$
and $\talg{\alpha}_1$ is a 1-form, then by definition of the Lie derivation we
have
$$
\lder{X_{\talg{\alpha}}}{\talg{\omega}}|_N=(d\varphi\cdot\talg{\omega}(X_{\talg{\alpha}_1}))|_N+(\varphi\cdot\lder{X_{\talg{\alpha}_1}}{\talg{\omega}})|_N
$$
The first term of the RH side is 0 because the vector field
$X_{\talg{\alpha}_1}|_N$ is always tangent to the symplectic leaf $N$, and
$\talg{\omega}|_N=\omega$ is a section of $\ortho$. Therefore we have
$$
\lder{X_{\talg{\alpha}_1}}{\talg{\omega}}|_N=\varphi|_N\cdot\lder{X_{\talg{\alpha}_1}}{\talg{\omega}}|_N
$$
from which easily follows as the independence of the expression
$\lder{X_{\talg{\alpha}}}{\talg{\omega}}|_N$ from the extension $\talg{\alpha}$
so the $\smooth{N}$-linearity of $\nabla_\alpha(\omega)$ by the argument
$\alpha$.

\newcommand{\zxa}{X_{\talg{\alpha}}}
For $\talg{\omega}=\varphi\cdot\talg{\omega}_1$, we have
\begin{equation}\label{leibniz}
\lder{\zxa}{\varphi\cdot\talg{\omega}_1}|_N=(\zxa(\varphi)\cdot\talg{\omega}_1)|_N+(\varphi\cdot\lder{\zxa}{\talg{\omega}_1})|_N
\end{equation}
From this we obtain the independence from the extension $\talg{\omega}$ in the
following way: if $\talg{\omega}|_N=\talg{\omega}'|_N$, then
$\talg{\omega}-\talg{\omega}'=\sum\varphi_i\talg{\omega}_i$, where
$\varphi_i|_N=0$; then recall that as $\zxa|_N$ is tangent to $N$, we have that
$\zxa(\varphi_i)|_N=0$. From the equality \ref{leibniz} also follows the
Leibniz rule for the argument $\omega$ in $\nabla_\alpha(\omega)$.

Let us describe the dual definition of linear connection, which is more
compatible with the algebraic definition introduced in the previous section. In
this case, instead of the vector bundle $\ortho$, we take the quotient bundle
$\frac{T(M)|_N}{T(N)}\equiv V(N)$, which is the bundle of dual spaces of the
fibers of $\ortho$. For a section $\alpha$ of $\restr$ and a section $s$ of
$V(N)$, we set
\begin{equation}\label{bott}
\nabla_\alpha(s)=\rho([\zxa,\talg{s}])
\end{equation}
where: $\talg{\alpha}$ is an extension of $\alpha$ as in the previous
definition; $\talg{s}$ is a vector field on $M$, such that
$\rho(\talg{s}|_N)=s$ and $\rho:T(M)|_N\To V(N)$ is the qoutient map. The
independence of the expression \ref{bott} from the extensions $\talg{\alpha}$
and $\talg{s}$, as the properties of covariant derivation easily follow from
the definition by analogy to the case of the dual definition.

In general, if $\nabla$ is a covariant derivation on some vector bundle, its
dual on the dual vector bundle is defined by the equality
\begin{equation}\label{dualconns}
\scalar{\nabla_\alpha(\omega)}{s}=\alpha\scalar{\omega}{s}-\scalar{\omega}{\nabla_\alpha(s)}
\end{equation}
where $\omega$ and $s$ are sections of the bundles dual to each other. Keeping
in mind this equality and follow the definitions of the linear connections on
$\ortho$ and its dual $V(N)$, we obtain
$$
\lder{\zxa}{\talg{\omega}}(\talg{s})=\talg{s}\pars{\talg{\omega}(\zxa)}+
\zxa\pars{\talg{\omega}(\talg{s})}-\talg{s}\pars{\talg{\omega}(\zxa)}-
\talg{\omega}([\zxa,\talg{s}])=\zxa(\talg{\omega}(\talg{s}))-\talg{\omega}([\zxa,\talg{s}])
$$
from which follows the equality \ref{dualconns} for $\omega$ -- a section of
$\ortho$, $s$ -- a section of $V(N)$ and $\alpha\in\restr$.

Let us consider an example which is rather simple but useful for the
demonstration of various interesting properties of singular Poisson structures.
\begin{exmp}\label{example}
Let $S$ be a 2-dimensional symplectic manifold, with Poisson bracket $\brck{\
}{\ }$. Let $\varphi$ be a smooth function on $S$, such that
$\varphi^{-1}(0)=\set{x_0}$. Conider a modified Poisson bracket $\brck{\ }{\
}_1=\varphi\cdot\brck{\ }{\ }$ (as the manifold is 2-dimensional this bracket
satisfies all the properies required for a Poisson bracket). The symplectic
foliation for the modified bracket consists of two symplectic leaves:
$\set{x_0}$ and $S\setminus\set{x_0}$. A fiber of the transversal foliation for
the leaf $S\setminus\set{x_0}$ is just $\set{0}$. Hence it is more interesting
to consider the linear connection for the leaf consisting of just one point
$\set{x_0}$. In this case the transversal space is $V(x_0)=T_{x_0}(S)$ and
thus, by definition, the covariant derivation, corresponding to the linear
connection, is an action of the elements of $T_{x_0}^*(S)$ on $T_{x_0}(S)$:
\begin{equation}\label{bott1}
\nabla_\alpha(V)=[X^1_f,\talg{V}]|_{x_0}=[\varphi\cdot
X_f,\talg{V}]|_{x_0}=V(\varphi)\cdot X_f(x_0)
\end{equation}
where $f\in\smooth{S}$ is such that $df(x_0)=\alpha$; $X^1_f$ is the
Hamiltonian vector field for $f$ with respect to the modified bracket $\brck{\
}{\ }_1$ and $X_f$ is the Hamiltonian vector field for $f$ with respect to the
original bracket $\brck{\ }{\ }$. Notice that as the original bracket is
nondegenerate, we have that $\alpha\neq0\tolfas X_\alpha\neq0$, therefore the
equation $\nabla_\alpha(V)=0$, for $V$, when $\alpha\neq0$ is equivalent to
$V\in\ker(\varphi'(x_0))$. In other words, the ``flat'' sections of the vector
bundle $V(x_0)$ over $\set{x_0}$ are the elements of $\ker(\varphi'(x_0))$.
Also it follows from the equality \ref{bott1} that the linear connection on the
leaf $\set{x_0}$ is flat if and only if $\varphi'(x_0)=0$.
\end{exmp}
%%%%%%%%%%%%%%%%%%%%%%%%%%%%%%%%%%%%%%%%%%%%%%%%%%%%%%%%%%%%%%%
\section{Linear Connection for Transversal Mulivectors and Schouten-Nijenhuis Bracket}\label{bott-for-trans-multi}
For an arbitrary smooth manifold M we denote by $\fields{k}{M},\
k=0,\ldots,\infty$, the space of skew-symmetric contravariant tensor
(multivector) fields of degree $k$ on $M$ ($\fields{0}{M}=\smooth{M}$). The
Schouten-Nijenhuis bracket is a linear operation
$$
[\ ,\ ]:\fields{p}{M}\times\fields{q}{M}\To\fields{p+q-1}{M}
$$
with the following properties (see \cite{Nijen})
\begin{equation}\label{schout-props}
\begin{array}{ll}
(1) & [U,V]=(-1)^{|U|\cdot|V|}\cdot[V,U] \\
(2) & [U,V\wedge W]=[U,V]\wedge W+(-1)^{(|U|+1)\cdot|V|}\cdot V\wedge[U,W] \\
(3) & (-1)^{|U|\cdot(|W|-1)}\cdot[U,[V,W]]+(-1)^{|V|\cdot(|U|-1)}\cdot[V,[W,U]]+ \\
    & +(-1)^{|W|\cdot(|V|-1)}\cdot[W,[U,V]]=0
\end{array}
\end{equation}
For monomial type multivector fields the formula for the Schoute-Nijenhuis
bracket is
\begin{equation}\label{schout}
\begin{array}{c}
[X_1\wedge\cdots\wedge X_m,Y_1\wedge\cdots\wedge Y_n]=\\
\\
(-1)^{m+1}\jami{i,j}{}(-1)^{i+j}[X_i,Y_j]\wedge
X_1\wedge\cdots\wedge\hat{X_i}\wedge\cdots\wedge X_m\wedge
Y_1\wedge\cdots\wedge\hat{Y_j}\wedge\cdots\wedge Y_n
\end{array}
\end{equation}
If we define a generalization of the ordinary Lie derivation for the
multivector fields as
$$
L_X=i_X\circ d-(-1)^{|X|}d\circ i_X
$$
then there is the following relation between this operation and the
Schouten-Nijenhuis bracket
\begin{equation}\label{lieschout}
[L_X,i_Y]=i_{[X,Y]}
\end{equation}
where $X$ and $Y$ are multivector fields and $i$ denotes the inner product of a
multivector field and a differential form.

For an arbitrary vector bundle $E$ and a covariant derivation $\nabla$ on it,
the canonical extension of $\nabla$ to the Grassmann algebra bundle
$\wedge(E)=\jami{k=0}{\infty}\wedge^k(E)$ is defined by the formula
\begin{equation}\label{covderext}
\nabla_X(S_1\wedge\cdots\wedge S_k)=\sum(-1)^{i+1}\nabla_X(S_i)\wedge
S_1\wedge\cdots\wedge\hat{S_i}\wedge\cdots\wedge S_k
\end{equation}
where $S_1\wedge\cdots\wedge S_k$ is a section of $\wedge^k(E)$.

For a symplectic manifold $M$ and a symplectic leaf $N\subset M$, the $k$-th
exterior degree of the vector bundle $V(N)=\frac{T(M)|_N}{T(N)}$ is canonically
isomorphic to $\frac{\wedge^kT(M)|_N}{T(N)\wedge(\wedge^{k-1}T(M))|_N}$. Here
$T(N)\wedge(\wedge^{k-1}T(M))|_N$ is the intersection of the ideal generated by
$T(N)$ in the Grassmann algebra bundle $\wedge T(M)|_N$, with
$\wedge^kT(M)|_N$. According to the formulas \ref{schout} and \ref{covderext},
the extension of the Bott covariant derivation on $V(N)$ to $\wedge^kV(N)$ can
be defined as
\begin{equation}\label{BottSchouten}
\nabla_\alpha(U)=\rho([\zxa,\talg{U}]|_N)
\end{equation}
where $\alpha$ is a section of $\restr$; $U$ is a section of $\wedge^kV(N)$;
$$
\rho:\wedge^kT(M)|_N\To\wedge^kV(N)
$$
is the quotient map; $\talg{\alpha}$ is an extension of $\alpha$; $\talg{U}$ is
such a section of $\wedge^kT(M)$ (i.e., a multivector field on $M$) that
$\rho(\talg{U}|_N)=U$ and $[\ ,\ ]$ is the \emph{Schouten-Nijenhuis Bracket}.
\begin{exmp}\label{example2}
Let $S$, $\brck{\ }{\ }_1=\varphi\brck{\ }{\ }$ and $x_0$ be the same objects
as in the Example \ref{example}. We have that
$\wedge^2V(x_0)=\wedge^2T_{x_0}(S)$. For $\alpha\in T^*_{x_0}(S)$ and $U\wedge
V\in T_{x_0}(S)$, we have
$$
\begin{array}{l}
\nabla_\alpha(U\wedge V)=\nabla_\alpha(U)\wedge V-\nabla_\alpha(V)\wedge U=\\
\\
U(\varphi)\cdot X_\alpha\wedge V-V(\varphi)\cdot X_\alpha\wedge
U=X_\alpha\wedge(U(\varphi)\cdot V-V(\varphi)\cdot U)
\end{array}
$$
After this suppose that we need to solve the equation $\nabla_\alpha(U\wedge
V)=0$ for $U\wedge V$. As $\dim(S)=2$, the space $\wedge^2T_{x_0}(S)$ is
1-dimensional. Therefore we are free to take the vector $V$ as an element of
$\ker(\varphi'(x_0))$. Hence we obtain
$$
\nabla_\alpha(U\wedge V)=U(\varphi)\cdot X_\alpha\wedge V=0
$$
From this follows that if $\varphi'(x_0)=0$, then the extended linear
connection on $\wedge^2T_{x_0}(S)$ is flat, otherwise, the equation
$\nabla_\alpha(U\wedge V)=0$ has a nontrivial solution if and only if
$X_\alpha\in\ker(\varphi'(x_0))$. The both case are summarized as the
following: the equation $\nabla_\alpha(U\wedge V)=0$ for has nontrivial
solution if and only if $X_\alpha\in\ker(\varphi'(x_0))$, and the set of
solutions is $U\wedge\ker(\varphi'(x_0))$ (of course
$\cong\wedge^2T_{x_0}(S)$).
\end{exmp}
%%%%%%%%%%%%%%%%%%%%%%%%%%%%%%%%%%%%%%%%%%%%%%%%%%%%%%%%%%%%%%%%%%%%%%%%%%%%%%%%%%%%%%%%%%%%
\section[Lie Algebra Cohomologies with Values in a Module]{\large{Lie Algebra Cohomologies with Values in a Module}}\label{section1}
We start from review of some definitions and facts from the cohomology theory
of Lie algebras. Let $L$ be a Lie algebra, $A$ be an associative and
commutative algebra over the field of complex (or real) numbers and $S$ be a
module over the algebra $A$.
\begin{defn}\label{Lie_module_definition}
The triple $(L,S,A)$ is said to be a Lie module over the Lie algebra $L$ if the
following conditions are satisfied:
\begin{itemize}
\item[\textbf{(l1)}]
$L$ is a module over the algebra $A$ and there is a Lie algebra homomorphism
$\phi:L\To Der(A)$, which is also a homomorphism of $A$-modules, such that for
any $a\in A$ and $X,Y\in L$, we have: $[X,aY]=X(a)Y+a[X,Y]$. Here $Der(A)$
denotes the space of derivations of the algebra $A$ and $X(a)$ denotes
$\phi(X)(a)$.
\item[\textbf{(l2)}]
There is a Lie algebra homomorphism $\psi:L\To End_\Complex(S)$, which is also
a homomorphism of $A$-modules, such that for any $X\in L,\;s\in S$ and $a\in
A$, we have: $X(as)=X(a)s+aX(s)$. Here $End_\Complex(S)$ denotes the space of
$\Complex$-linear mappings from $S$ to itself and $X(s)$ denotes $\psi(X)(s)$.
\end{itemize}
\end{defn}
Sometimes, for brevity, the term ``Lie module over $L$'' will be used for the
$A$-module $S$.

For any integer number $m\geq1$, let us denote the space of $A$-multilinear
mappings from $\underbrace{L\times\cdots\times L}_{m-times}$ to the Lie module
$S$, by $\Omega_A^m(L,S)$. We also set that $\Omega_A^0(L,S)=S$ and
$\Omega_A(L,S)=\bigoplus\limits_{m=0}^{\infty}\Omega_A^m(L,S)$.

For any $\omega\in\Omega_A^m(L,S)$ and $\set{X_1,\cdots,X_{m+1}}\subset L$,
define $d\omega$ by the well-known Koszul formula:
$$
\begin{array}{l}
(d\omega)(X_1,\cdots,X_{m+1})=\jami{i}{}(-1)^{i-1}X_i\omega(X_1,\cdots,\widehat{X_i},\cdots,X_{m+1})+\\
+\jami{i<j}{}(-1)^{i+j}\omega([X_i,X_j],\cdots,\widehat{X_i},\cdots,\widehat{X_j},\cdots)
\end{array}
$$
\begin{lem}\label{form_differential_is_form_lemma}
If the $A$-module $L$ is projective then the following three conditions are
equivalent:
\begin{enumerate}
\item
For any $a\in A$ and $s\in\Omega^0_A(L,S)=S$: $d(as)=d(a)s+ad(s)$.
\item
The $A$-module $S$ is a Lie module over $L$.
\item
For any integer $m\geq1$ and $\omega\in\Omega^m_A(L,S)$, the mapping
$$d\omega:L^{m+1}\To S$$
is $A$-multilinear (i.e.,
$(d\omega)(aX_1,\cdots,X_{m+1})=a(d\omega)(X_1,\cdots,X_{m+1})$).
\end{enumerate}
\end{lem}
\begin{prf}
It is clear that the conditions 1 and 2 are equivalent, because for any $X\in
L$, we have that:
$$
\pars{d(as)}(X)=X(as)=X(a)s+aX(s)=\pars{d(a)s+ad(s)}(X).
$$
Also, it can be verified by direct calculation that from the condition 2 (or,
which is the same, from 1) follows the condition 3. Now, suppose that the
condition 3 is true. In this case For any $X,Y\in L,\;a\in A$ and
$\omega\in\Omega_A^1(L,S)$, we have:
$$
\begin{array}{l}
(d\omega)(X,aY)=a(d\omega)(X,Y)\quad\Rightarrow\\
\\
X\omega(aY)-aY\omega(X)-\omega([X,aY])=aX\omega(Y)-aY\omega(X)-a\omega([X,Y])\;\Rightarrow\\
\\
X(a\omega(Y))-aY\omega(X)-X(a)\omega(Y)-a\omega([X,Y])=\\
\\
=aX\omega(Y)-aY\omega(X)-a\omega([X,Y])\qquad\Rightarrow\\
\\
X(a\omega(Y))=X(a)\omega(Y)+aX(\omega(Y))
\end{array}
$$
Because it is assumed that the $A$-module $L$ is projective, for any $s\in S$
can be found such $\omega\in\Omega^1_A(L,S)$ and $Y\in L$, that $\omega(Y)=s$.
Therefore, we obtain that $X(as)=X(a)s+aX(s)$, i.e., $S$ is a Lie module over
the Lie algebra $L$.
\end{prf}

Hence, we have that if the $A$-module $S$ is a Lie module over the Lie algebra
$L$, the operator $d$, carries the space $\Omega^m_A(L,S)$ into
$\Omega^{m+1}_A(L,S)$, which implies that the pair $(\Omega_A(L,S),d)$ is a
differential complex.
%%%%%%%%%%%%%%%%%%%%%%%%%%%%%%%%%%%%%%%%%%%%%%%%%%%%%%%%%%%%%%%%%%%%%%%%%%%%%%%%%%%%%%%
\section[Lie Algebra Homologies and Supercommutator]{\large{Lie Algebra Homologies and Supercommutator}}\label{section2}
\newcommand{\g}{\ensuremath{\mathfrak{g}} }
Let \g be a Lie algebra and a module over a commutative and associative algebra
$B$. If the Lie algebra bracket in \g is bilinear for the elements of $B$
($[x,ay]=a[x,y],\;\all{x,y}\in\g$ and $a\in B$) one has a well-defined boundary
operator $\delta:\wedge_B^m(\g)\To\wedge_B^{m-1}(\g)$, where $\wedge_B$ denotes
the exterior product of a $B$-module by itself:
$$
\begin{array}{l}
\delta(x_1\wedge\cdots\wedge x_m)=\jami{i<j}{}(-1)^{i+j}[x_i,x_j]\wedge\cdots\wedge\widehat{x_i}\wedge\cdots\wedge\widehat{x_j}\wedge\cdots\wedge x_m,\;m>1\\
\textrm{ and }\quad\delta(x)=0,\;\all{\set{x,x_1,\cdots,x_m}}\subset\g
\end{array}
$$
The homology given by the boundary operator $\delta$ is known as the homology
of the Lie algebra \g with coefficients from the algebra $B$.

There is a useful relation between the coboundary operator $\delta$ and the
Schouten-Hijenhuis bracket on the exterior algebra $\wedge_B(\g)$:
$\all{u}\in\wedge_B^m(\g)$ and $\all{v}\in\wedge_B(\g)$, we have
\begin{equation}\label{super_as_oper_deviation_formula}
[u,v]=\delta(u)\wedge v+(-1)^mu\wedge\delta(v)-\delta(u\wedge v)
\end{equation}
It can be said that the supercommutator $[\cdot,\cdot]$ measures the deviation
of the boundary operator $\delta$ from being an antidifferential of degree -1
(see \cite{Me}). From the formula \ref{super_as_oper_deviation_formula}, easily
follows that the induced supercommutator on the homology space
$H_{\bullet}(\g,B)$ is trivial: for $u\in\wedge_B^m\g$ and $v\in\wedge_B\g$,
such that $\delta(u)=\delta(v)=0$, we have $[u,v]=-\delta(u\wedge v)$, which
implies that the homology class of the element $[u,v]$ is trivial.

From the formula \ref{super_as_oper_deviation_formula}, also follows that the
homology space $H_\bullet(\g,B)$ does not inherits the exterior algebra
structure from $\wedge_B(\g)$: if the elements $u,v\in\wedge_B(\g)$ are closed,
then we have that $\delta(u\wedge v)=-[u,v]$.
%%%%%%%%%%%%%%%%%%%%%%%%%%%%%%%%%%%%%%%%%%%%%%%%%%%%%%%%%%%%%%%%%%%%%%%%%%%%%%%%%%%%%%%
\section[The Homology Space of a Lie Algebra Ideal]{\large{The Homology Space of a Lie Algebra Ideal}}\label{section3}
Let $L$ be a Lie algebra and $A$ be a commutative and associative algebra over
$\Complex$ (or $\Real$). Suppose that the pair $(L,A)$ satisfies the condition
\emph{l1} from the definition \ref{Lie_module_definition}.

Let $V\subset L$ be an ideal in the Lie algebra $L$ (i.e., $\all{v}\in V$ and
$\all{x}\in L$: $[v,x]\in V$). Define the subalgebra $A_V$, in the algebra $A$,
as
$$
A_V=\set{a\in A\,|\,v(a)=0,\;\all{v}\in V}
$$
and the subspace $V'$ in $L$, as
$$
V'=\set{\jami{i=1}{n}a_iv_i\,\vbar\,n\in\mathbb{N}\;a_i\in A_V,\;v_i\in
V,\;i=1,\cdots,n}
$$
The latter can be defined as the minimal $A_V$-submodule of $L$, containing
$V$.

For any $x\in L,\;v\in V$ and $a\in A_V$, we have
$$
v(x(a))=[v,x](a)+x(v(a))=0\quad\imply\quad x(a)\in A_V
$$
which implies that the subalgebra $A_V$ is invariant under the action of the
elements of $L$. Therefore, we have that $[x,av]=x(a)v+a[x,v]\in V'$, which
means that $V'$ is also an ideal in the Lie algebra $L$. We call the ideal $V'$
the \emph{complement} of the ideal $V$.

From the fact that $V$ is a subspace of $V'$ follows that $A_{V'}\subset A_V$,
but on the other hand, for any $a,b\in A_V$ and $v\in V$, we have that
$(bv)(a)=b\cdot v(a)=0$, which means that any $a\in A_V$ is also an element of
$A_{V'}$. Hence, the two algebras $A_V$ and $A_{V'}$, coincide.

The ideal $V$ in the Lie algebra $L$ will be called \emph{complete}, if $V=V'$.
Further, by default, we assume that the ideal $V$ is complete (i.e., $V$ is a
module over the the algebra $A_V$).

Consider the homology space of the Lie algebra $V$ with coefficients from
$A_V$. One can define an action of the Lie algebra $L$ on the $A_V$-module
$H_\bullet(V,A_V)$: for $X\in L$ and $v\in V$ such that $\delta(v)=0$, let
$X([v])=[[X,v]]$, where $[\cdot]$ denotes the homology class and
$[\cdot,\cdot]$ denotes the supercommutator.
\begin{prop}
The action of the Lie algebra $L$ on the homology space\\
$H_\bullet(V,A_V)$ is correctly defined and gives a Lie module structure on the\\
$A_V$-module $H_\bullet(V,A_V)$.
\end{prop}
\begin{prf}
As it follows from the formula \ref{super_as_oper_deviation_formula}, we have
that for any $X\in L$ and a closed element $u\in\wedge_{A_V}(V)$:
$[X,u]=-\delta(X\wedge u)$, which implies that $[X,u]$ is closed. If the
element $u$ is exact and $u=\delta(v)$, then we have
$$
\delta([X,v])=\delta(-X\wedge\delta(v)-\delta(X\wedge v))=-\delta(X\wedge
u)=-[X,u]
$$
which implies that the element $[X,u]$ is exact. From these two fact follows
that the closed and exact element in the exterior algebra $\wedge_{A_V}(V)$ are
invariant under the supercommutator with the elements of the Lie algebra $L$.
Therefore, the action, $X([u])=[[X,u]]$, of $L$ on the homology space
$H_\bullet(V,A_V)$ is correctly defined. By definition of the algebra $A_V$ and
the properties of the supercommutator, easily follow that for any $a\in A_V$:
$(aX)([u])=aX([u])$ and $X(a[u])=X([au])=[[X,au]]=X(a)[u]+aX([u])$; which
implies that the $A_V$-module $H_\bullet(V,A_V)$ is a Lie module over the Lie
algebra $L$.
\end{prf}\bigskip\\
\newcommand{\homid}[3][]{\Omega_{A_V}^{#1}(#2,\;H_{#3}(V,A_V))}
\newcommand{\hmid}[3][]{\Omega_{A_V}^{#1}(#2,H_{#3}(V,A_V))}
It follows from the above proposition, that for each integer $n\geq0$, one can
consider the the following differential complex $\pars{\homid{L}{n},\;d}$,
which gives the Lie algebra cohomology of $L$ with values in the $A_V$-module
$H_n(V,A_V)$.

As it was mentioned early, the induced supercommutator on each homology space
$H_n(V,A_V),\;n=0,\cdots,\infty$, is trivial. Therefore, the action of the Lie
algebra $V$ on each $H_n(V,A_V)$ is trivial, which implies that the submodule
of the $A_V$-module $\homid{L}{\bullet}$, consisting of such forms $\omega$,
that $i_v(\omega)=0,\;\all{v}\in V$, is invariant under the action of the
differential $d$ ($i_v\omega=0\;\imply\;i_v(d\omega)=0,\;\all{v}\in V$). Hence,
one can consider the subcomplex $\pars{\homid{L}{\bullet}_0,\;d}$ of the
differential complex $\pars{\homid{L}{n},\;d}$, where
$$
\homid{L}{\bullet}_0=\set{\omega\in\homid{L}{\bullet}\;|\;i_v(\omega)=0,\;\all{v}\in
V}
$$
Actually, the complex $\pars{\homid{L}{\bullet}_0,\;d}$ is canonically
isomorphic to the complex $\pars{\homid{L/V}{\bullet},\;d}$.
%%%%%%%%%%%%%%%%%%%%%%%%%%%%%%%%%%%%%%%%%%%%%%%%%%%%%%%%%%%%%%%%%%%%%%
\section[The Characteristic Class of a Lie Algebra Ideal]{\large{The Characteristic Class of a Lie Algebra Ideal}}
Let us denote by $\homid{L}{\bullet}_1$ the quotient module
$$\homid{L}{\bullet}/\homid{L}{\bullet}_0$$
For any integer $n\geq 0$, we have the following short exact sequence:
$$
\begin{array}{ll}
0\;\To\;\homid{L/V}{n}\;\To\;\homid{L}{n} & \stackrel{\pi}{\To} \\
 \stackrel{\pi}{\To}\;\homid{L}{n}_1\;\To\;0 &
\end{array}
$$
which induces the standard homomorphism of cohomology spaces
$$
\sigma_n:H^\bullet(L,\;H_n(V,A_V))_1\To H^{\bullet\,+1}(L/V,\;H_n(V,A_V))
$$
where $H^\bullet(L,\;H_n(V,A_V))_1$ denotes the cohomology space of the complex
$(\homid{L}{n},\;d)$ and $\sigma_n$ is defined as:
$\sigma_n([\pi(\omega)])=[d\omega]$, for $\omega\in\homid{L}{n}$, such that
$d(\pi(\omega))=0$
($\;\imply\;\pi(d\omega)=0\;\imply\;d\omega\in\homid{L/V}{n}$); here
$[d\omega]$ denotes the cohomology class of the form $d\omega$ in
$H^{\bullet\,+1}(L/V,\;H_n(V,A_V))$.

\newcommand{\cohom}[2][]{\Omega_{A_V}^{#1}(#2,\,H_V)}

Consider the special case when $n=1$. Let us denote the homology space
$H_1(V,A_V)$ by $H_V$. So, we have an exact sequence
\begin{equation}\label{cohomology_exact_seq2}
0\;\To\;\cohom[\bullet]{L/V}\;\To\;\cohom[\bullet]{L}\;\stackrel{\pi}{\To}\;\cohom[\bullet]{L}_1\;\To\;0
\end{equation}
and the homomorphism: $\sigma_1:H^{\bullet}(L,H_V)_1\To
H^{\bullet\,+1}(L/V,H_V)$.

There is a one-to-one correspondence between the set of splittings of the
following short exact sequence
$$
0\;\To\;V\;\To\;L\;\stackrel{p}{\To}\;L/V\;\To\;0
$$
and the set of such homomorphisms $\alpha:L\To V$, that
$\alpha(v)=v,\;\all{v}\in V$ (projection operator). Any such projection
operator $\alpha$, defines a 1-form $\widetilde{\alpha}\in\cohom[1]{L}$:
$$
\widetilde{\alpha}(X)=[\alpha(X)],\;\all{X}\in L
$$
\begin{lem}
For any two projections $\alpha,\beta:L\To V$, the forms
$\pi(\widetilde{\alpha})$ and $\pi(\widetilde{\beta})$ in $\cohom[1]{L}_1$, are
equal and closed.
\end{lem}
\begin{prf}
For such $\alpha$ and $\beta$ we have: $(\alpha-\beta)(v)=v-v=0,\;\all{v}\in
V$, which implies that
$i_v(\widetilde{\alpha}-\widetilde{\beta})=0\;\imply\;\widetilde{\alpha}-\widetilde{\beta}\in\cohom[1]{L/V}\;
\imply\;\pi(\widetilde{\alpha})=\pi(\widetilde{\beta})$.\\
For any $v\in V$ and $X\in L$, we have:
$$
\begin{array}{l}
(d\widetilde{\alpha})(v,X)=v(\widetilde{\alpha}(X))-X(\widetilde{\alpha}(v))-\widetilde{\alpha}([v,X])=\\
\\
=-X([v])-[[v,X]]=-[[X,v]]-[[v,X]]=0
\end{array}
$$
Therefore, we obtain that $i_v(d\widetilde{\alpha})=0,\;\all{v}\in V$, which
means that the form $d\widetilde{\alpha}$ belongs to the submodule
$\cohom[2]{L/V}$, or equivalently:
$\pi(d\widetilde{\alpha})=d\pi(\widetilde{\alpha})=0\;\imply\;\pi(\widetilde{\alpha})$
is closed in $\cohom[1]{L}_1$.
\begin{rem}
$d\widetilde{\alpha}-d\widetilde{\beta}=d(\widetilde{\alpha}-\widetilde{\beta})\;\imply$
the forms $\widetilde{\alpha}$ and $\widetilde{\beta}$ are cohomological in
$\cohom[2]{L/V}$, and their cohomology class is exactly
$\sigma_1(\pi(\widetilde{\alpha}))\in H^2(L/V,\;H_1(V,A_V))$.
\end{rem}
\end{prf}\medskip\\
As it follows from the above lemma, for any projection homomorphism
$\alpha~:~L~\To~V$, the element $\pi(\widetilde{\alpha})$ in $\cohom[1]{L}_1$
is one and the same. Let us denote this element (and also its cohomology class
in $H^1(L,H_V)_1$) by $\Delta$. The element $\sigma_1(\Delta)$ in
$H^2(L/V,\;H_1(V,A_V))$ we call the \emph{characteristic class} of the Lie
algebra ideal $V$ in $L$.
\begin{rem}
It is clear that if V is an ideal in the Lie algebra L, then the space
$[V,V]=\set{\jami{i}{}a_i[u_i,v_i]\;\vbar\;a_i\in A_V,\,u_i,v_i\in V}$ is also
an ideal in $L$ (an in $V$, too). We can consider the Lie algebra
$\widetilde{L}=L/[V,V]$ and its commutative subalgebra $\widetilde{V}=V/[V,V]$,
which is also an ideal in $\widetilde{L}$. It is clear that the homology space
$H_1(V,A_V)$ is the same as $\widetilde{V}$ and the quotient
$\widetilde{L}/\widetilde{V}$ is the same as $L/V$. One can consider the
following short exact sequence
\begin{equation}\label{seq2}
0\;\To\;\widetilde{V}\;\To\;\widetilde{L}\;\To\;L/V\;\To\;0
\end{equation}
As the Lie algebra $\widetilde{V}$ is commutative, in this situation, for any
connection form $\alpha:\widetilde{L}\To\widetilde{V}$ (a splitting of the
short exact sequence \ref{seq2}), its curvature form coincides with
$d\alpha\in\Omega^2_{A_V}(L/V,\widetilde{V})$, and the cohomology class of the
form $d\alpha$ in $d\alpha\in\Omega^2_{A_V}(L/V,\widetilde{V})$ coincides with
the characteristic class of the ideal $V$. Therefore, to study the properties
of the characteristic class of a Lie algebra ideal, we can consider the case
when the ideal is commutative.
\end{rem}
% ----------------------------------------------------------------

\end{document}